\documentclass[twoside]{article}

\usepackage{amsmath,amsthm,amssymb}

\usepackage{newtxtext,newtxmath}

\usepackage[text={12.5cm,19cm},centering,paperwidth=17cm,paperheight=24cm]{geometry}

\usepackage[fontsize=10.4pt]{scrextend}

\usepackage{graphicx}
\def\titlerunning#1{\gdef\titrun{#1}}

\makeatletter
\def\author#1{\gdef\autrun{\def\and{\unskip, }#1}\gdef\@author{#1}}
\def\address#1{{\def\and{\\\hspace*{15.6pt}}\renewcommand{\thefootnote}{}\footnote{#1}}\markboth{\autrun}{\titrun}}
\makeatother

\def\email#1{email: \href{mailto:#1}{#1} }
\def\subjclass#1{\par\bigskip\noindent\textbf{Mathematics Subject Classification 2020.} #1}
\def\keywords#1{\par\smallskip\noindent\textbf{Keywords.} #1}

\def\Tr{{\rm{Tr}}}
\def\bC{{\mathbb C}}

\newtheorem{thm}{Theorem}[section]

\newtheorem{lem}[thm]{Lemma}%
 \numberwithin{equation}{section}

\usepackage[hyperfootnotes=false,colorlinks=true,allcolors=blue]{hyperref}

\def\tr{{\rm Tr}}

\def\ra{{\rightarrow}}
\def\bX{{\bf X}}
\def\bB{{\bf B}}
\def\bA{{\bf A}}
\def\bG{{\bf G}}
\newcommand{\E}{\mathbb E}

\newcommand{\Pp}{\mathbb P}

\def\bR{{\mathbb R}}

\begin{document}

\titlerunning{Bernoulli Random Matrices}

\title{\textbf{Bernoulli Random Matrices}}

\author{Alice~Guionnet}

\date{}

\maketitle

\address{Alice Guionnet : Universit\'e de Lyon, ENSL, CNRS,  France; \email{alice.guionnet@ens-lyon.fr}}

\begin{abstract}
Random Matrix theory has become a field on its own with a breadth of new results, techniques, and ideas  in the last thirty years. In these proceedings of the 8ECM 2021, I illustrate some of these advances by describing what is known about the spectrum and the eigenvectors of Bernoulli matrices.
\end{abstract}

\section{Introduction}
Jacques (or Jakob) Bernoulli (1654-1705) was a renowned Swiss mathematician who made  important contributions to probability theory and partial differential equations. He was the first to discover the number $e$. But his most famous result is, at least for probabilists, the first proof of the law of large numbers. To this end, he analyzed the concept of Bernoulli law, which is the simplest non trivial distribution you can think about as the sum of two Dirac masses. It is the distribution of 
a random variable $b$ which can only take 2 values $0$ and $1$. We denote 
$$p=\Pp(b=1)=1-\Pp(b=0)$$
A very common example is a coin that, once thrown, falls  either on head (modeled by the state $1$) or tail (modeled by $0$). Even if one would expect in general the probability of each event to be equal to $1/2$, it may well be rather $p\in (0,1)$ if the coin is rigged. In Ars Conjectandi, Bernoulli showed that if one throws such a coin independently a number $n$ of times, then, with large probability,  one should see approximately $pn$ heads if $n$ is large enough. To state this law of large number more precisely, he showed that if $b_1,\ldots,b_n$ denotes the outcome of $n$ independent Bernoulli trials, then for any $a<p<b$
$$\lim_{n\rightarrow \infty} \Pp\left(\frac{1}{n}\sum_{i=1}^n b_i\in [a,b]\right)=1\,.$$
But how much can we choose $a,b$ close to $p$ so that this result remains true? Few years later, A. de Moivre (1667-1754) quantified the size of the error and proved the first central limit theorem, namely that
$a,b$ can be at  a distance of about $1/\sqrt{n}$ of $p$ in the sense that
$$\lim_{n\rightarrow \infty} \Pp\left(\frac{1}{\sqrt{n p(1-p)}}\sum_{i=1}^n (b_i-p)\in [a,b]\right)=\frac{1}{\sqrt{2\pi}} \int_a^b e^{-\frac{x^2}{2} }dx\,.$$
This was the first occurrence of the central limit theorem and the start of modern probability theory and statistics. Implicitly, we so far assumed that $p$ does not depend on $n$ and belongs to $(0,1)$. Later on, we shall also be interested in the case where $p$ depends on $n$. Then, it can be checked that the central limit theorem still holds as long as $pn$ goes to infinity. If $pn$ goes to a finite constant $c$, then it can not hold since $\sum_{i=1}^n b_{i}$ is an integer so that the above random variable is discrete. In fact, it converges towards the Poisson distribution
$$\lim_{n\rightarrow \infty} \Pp\left(\frac{1}{\sqrt{n p(1-p)}}\sum_{i=1}^n (b_i-p)\in [a,b]\right)=\sum_{k\in c+\sqrt{c}[a,b]}\frac{1}{k!} c^k e^{-c}$$
We will see later that this transition between such continuous and discrete limits is also key to describing the spectrum of Bernoulli Random Matrices. 
The last concept which is central in probability theory and important in these notes is entropy. It was introduced by Ludwig Boltzmann (1844--1906) and Claude Shannon (1916--2001) in physics and information theory respectively as a way to measure disorder.  For again $n$  independent Bernoulli trials with parameter $p$, it is defined  for any $q\in [0,1]$ by 

$${{\lim_{\varepsilon\downarrow 0}\lim_{n\ra\infty}\frac{1}{n}\ln \mathbb P\left(\frac{1}{n} \sum_{i=1}^{n} b_{i}\in [ q-\varepsilon, q+\varepsilon]\right)= -S_{p}(q)}}$$
where $S_{p}(q)= \frac{q}{p}\ln \frac{q}{p}+\frac{1-q}{1-p}\ln \frac{1-q}{1-p}$ is the entropy or rate function. 

\medskip

In this survey, I will discuss Bernoulli Random Matrices. A Bernoulli Random Matrix is a $n\times n$ symmetric matrix with independent Bernoulli entries (modulo the symmetry constraint) whose size $n$ is going to infinity and will discuss the law of large numbers, the fluctuations, and the entropy for their spectrum and eigenvectors.  There are many motivations to study random matrices. The first goes back to Wishart who considered random matrices to study correlations in large data sets. Such questions are very modern, with the need to analyze larger and larger data sets and machine learning.  The second comes from physics and the work of Wigner and Dyson. They proposed to model the Hamiltonian of excited nuclei by random matrices, an idea which turned out to be quite successful as indeed real nuclei turned out to have  energy levels distributed like the eigenvalues of random matrices. But Bernoulli matrices are special among all other random matrices because they describe the adjacency matrix of a Erd\"os-R\' enyi graph $G(n,p)$. Indeed, the later is just a graph build on $n$ (labeled) vertices, with an edge drawn independently between each couple of vertices with probability $p$. Studying the eigenvalues of the adjacency matrix of a graph gives valuable geometric  information, such as the size of its boundary (expanders) or the number of specific configurations, such as  triangles, that it contains. One can also be interested in the combinatorial properties of such matrices and for instance focus on the probability that the matrix is singular, see e.g \cite{VaVu}. My viewpoint will be to investigate the properties of the  eigenvalues and eigenvectors of Bernoulli random matrices, as a particularly nice and well-documented example of random matrices. 

To simplify, I will restrict myself to symmetric Bernoulli matrices $\bB_n$ throughout these notes :
$$\bB_n(i,j)= \bB_{n}(j,i)$$
and assume that $(B_{n}(i,j), i\le j)$ follows a Bernoulli law with parameter $p$. Also, I will take $\bB_n(i,i)$ random, but could take it equal to zero without changing much the statements of most of the  results.

My goal is to understand the spectrum of $\bB_n$ as well as the properties of its eigenvectors as $n$ goes to infinity. One can easily guess that these properties should depend on the parameter $p$. Indeed, thinking about the Erd\"os-R\' enyi graph, one sees that the average degree of a vertex is $pn$. The graph will be very dense if $pn$ goes to infinity fast enough but sparse if it is finite.

\begin{center}\includegraphics[width=3.5cm]{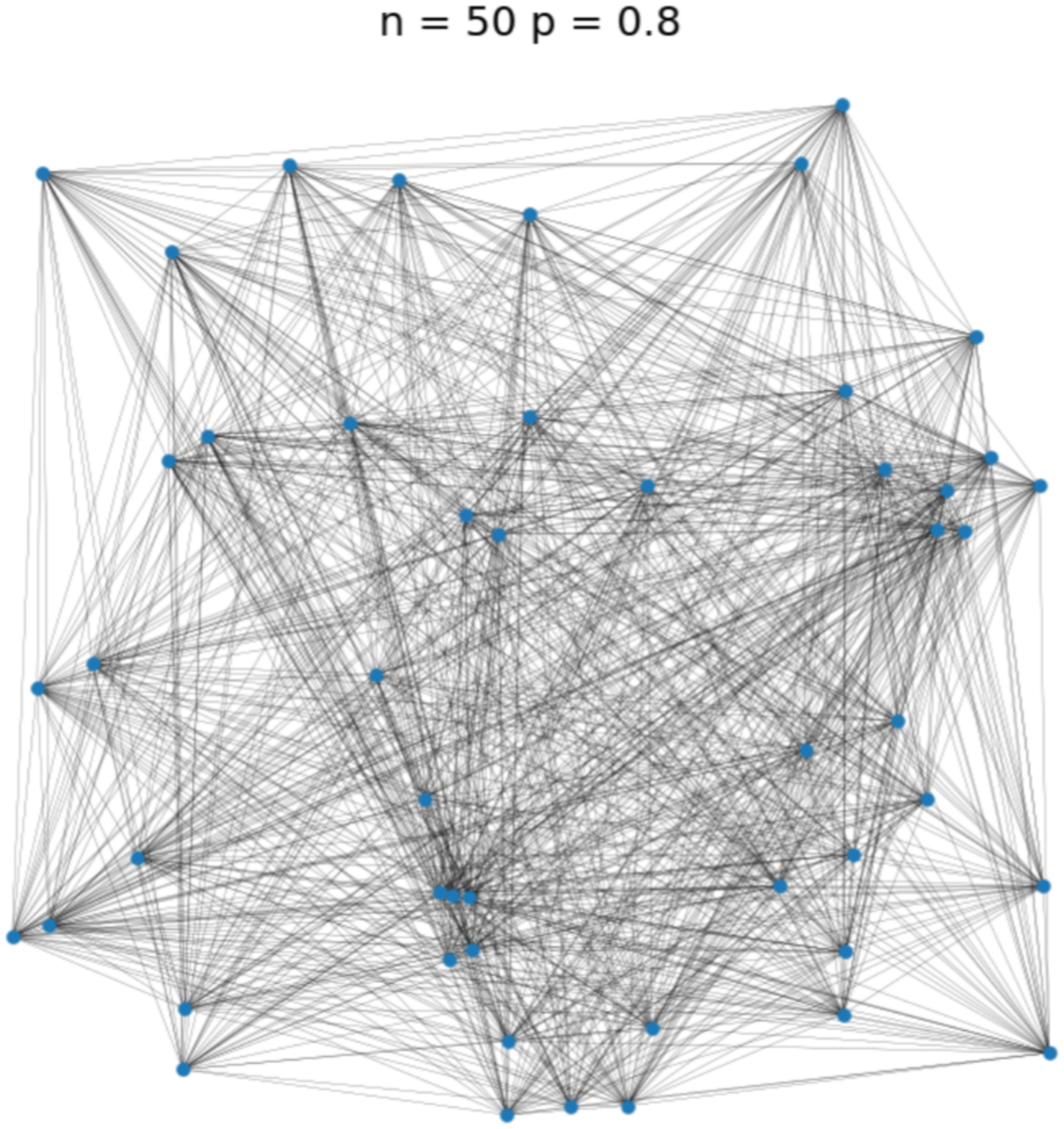}\hspace{0.5cm}\includegraphics[width=3.5cm]{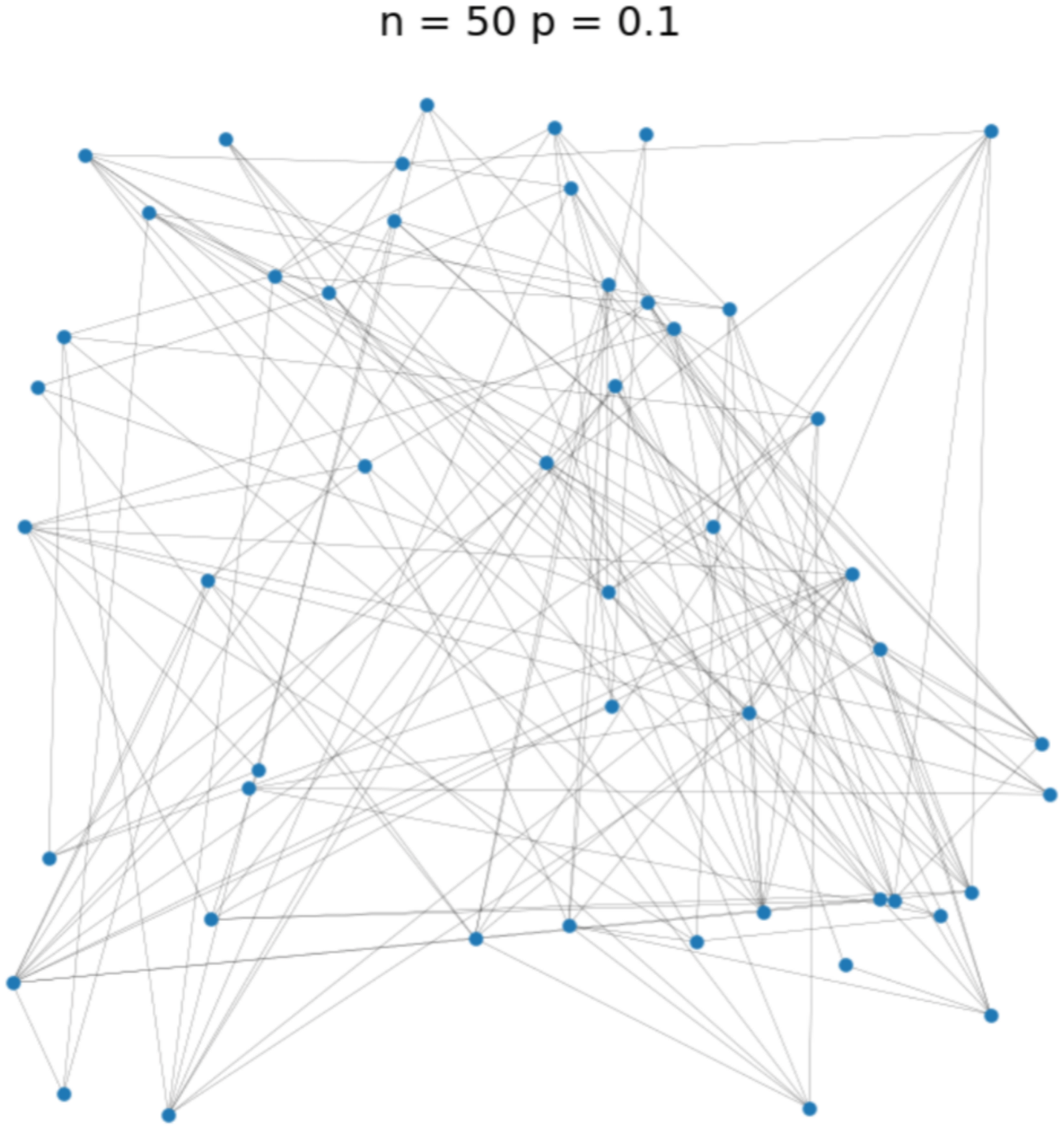}\hspace{0.5cm}\includegraphics[width=3.5cm]{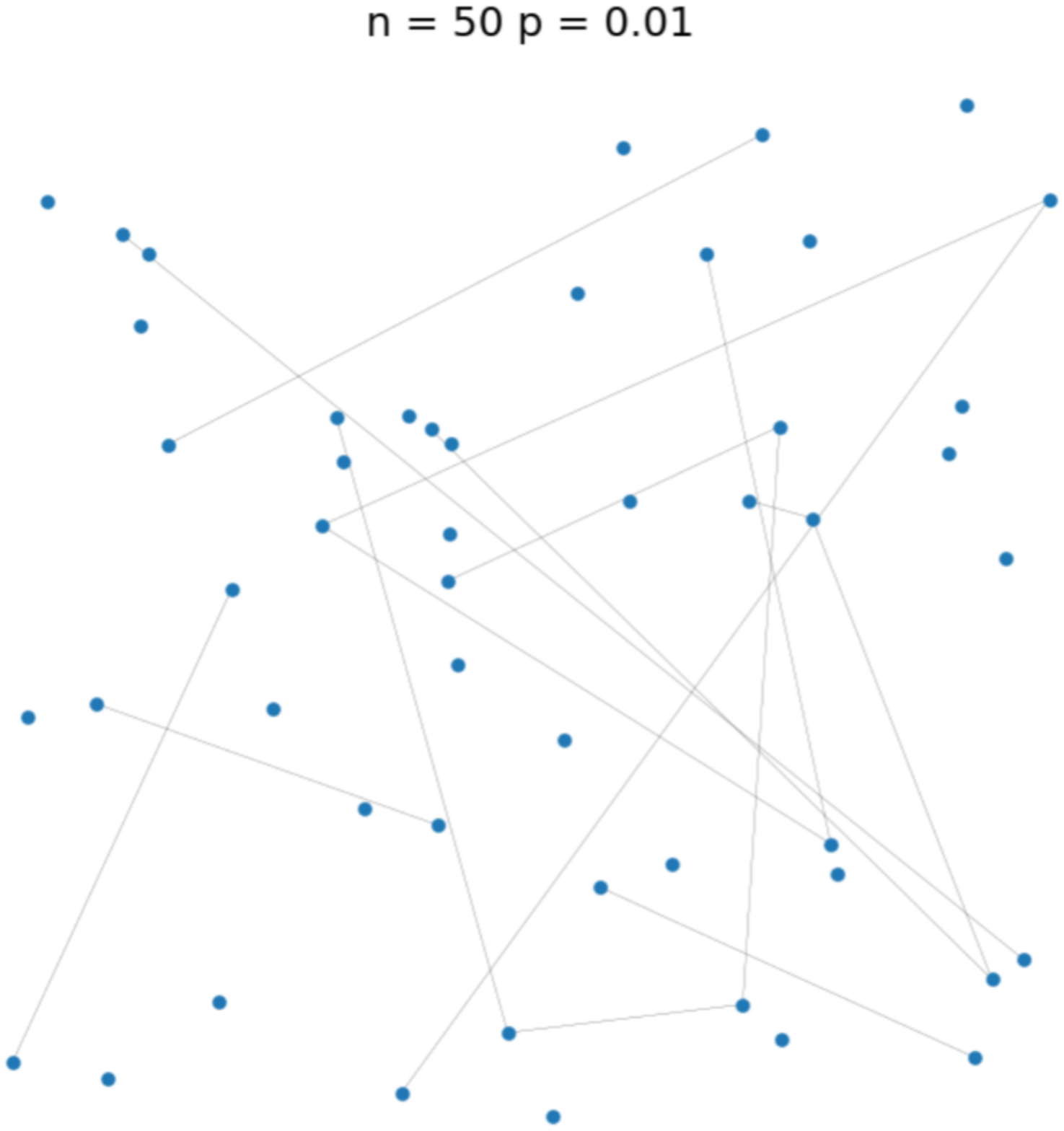}\end{center}
\begin{center}{\tiny Courtesy of D. Coulette}\end{center}

Indeed, it is well known  since the breakthrough paper of Erd\"os and R\'enyi that if $np<1$, $G(n,p)$ will almost surely have no connected component of size greater than $O(\log n)$, if $np=1$ there is a giant connected component but it is of size of order $n^{2/3}$, if $np$ goes to a constant $c>1$ it will have a unique giant component but lots of small components, and  isolated vertices will continue to exist until $ np>(1-\epsilon)\ln n$,
whereas if $np>(1+\epsilon)\ln n$ the graph will almost surely be connected. Here $\epsilon$ is some positive real number as small as wished. In the case where $np$ is of order $c$, the finite size connected components will create small diagonal blocks in the Bernoulli matrix, with entries equal either to zero or one and therefore finitely many possible eigenvalues. Hence, we expect the spectrum to accumulate at these possible values. But should there be other possible eigenvalues? Similarly, we see that the eigenvectors related with these eigenvalues are localized on a  few vertices. But should we also have delocalized eigenvectors ?  On the contrary, in the case where $np>(1+\epsilon)\ln n$, we may expect eigenvectors to be delocalized and the spectrum to be nicely continuous. In this case, a whole theory has been developed to show that the spectrum and the eigenvalues of Bernoulli matrices have the same properties as those of a random matrix with Gaussian entries. The latter is well known to be much easier to study, for instance, because the joint law of its eigenvalues is rather simple and independent of the eigenvectors.  
Conversely, Bernoulli matrices resemble more heavy tails matrices when $pn$ is of order one, in the sense that  it has mostly very tiny entries but a  few large  entries. Understanding the transition between these two behaviors is at the heart of random matrix theory. 

In this survey, I will start discussing the asymptotic behavior of the spectrum in both sparse and dense cases. Then, I will consider its fluctuations, both local and global, as well as the properties of its eigenvectors. Finally, I will discuss the large deviations of the spectrum, for instance how to estimate the probability that the second eigenvalue of Bernoulli matrices takes an unexpected value. 

{\bf Acknowledgements:} I thank Charles Bordenave, Christophe Garban, Jiaoyang Huang, Antti Knowles, Justin Salez and Ofer Zeitouni for pointing out recent developments of the subject and very helpful comments on preliminary versions of these proceedings. {This work was supported in part by ERC-2019-ADG Project 884584-LDRAM}.
\section{Law of large numbers}
In this section we shall see that the limiting distribution of the spectrum differs a lot according to the whether  $pn$ goes to infinity or not.

A first remark should be made about the matrix $\bB_n$: its entries are not centered. It will be more convenient to center them and renormalize the matrix properly. To this end we make the following decomposition
$$\bB_n= \sqrt{np(1-p)} \bX_n+ p {\bf 1}$$
where ${\bf 1}$ is a matrix whose entries are all equal to one, whereas the entries of $\bX_n$ are centered and renormalized to have covariance $1/n$:
$$\bX_n(ij)=\frac{\bB_n(ij)-p}{\sqrt{n p(1-p)}}\,.$$
The matrix ${\bf 1}$ has one non trivial eigenvalue which equals $n$, and flat eigenvector ${\bf 1}=(1/\sqrt{n},1/\sqrt{n},\ldots, 1/\sqrt{n})$. 
Conversely the spectrum of  $\bX_n$ has eigenvalues mostly of order one in the sense that  $\mathbb E[\Tr(\bX_n^2)]=\mathbb E[\sum \lambda_i^2]=n$.
Therefore the above decomposition shows that $\bB_n$ has a very large eigenvalue of order $n$, and the rest is given by the eigenvalues of $\bX_n$ taken on $ {\bf 1}^\perp$. Moreover, by Weyl's interlacing properties, the eigenvalues $(\lambda_i^B)_{1\le i\le n}$ of $\bB_n/ \sqrt{np(1-p)}$ and $(\lambda_i^X)_{1\le i\le n}$ of $  \bX_n$ are interlaced :
$$\lambda_n^X\le \lambda_n^B\le\lambda_{n-1}^X\cdots\le \lambda_1^X\le \lambda_1^B$$
Therefore, it is in general not difficult to retrieve the properties of the eigenvalues of $\bB_n/ \sqrt{np(1-p)}$ from those of $\bX_n$. Hereafter we will therefore concentrate mostly on $\bX_n$.
\subsection{Dense case}
The first result describes the asymptotic distribution of the spectrum in the dense case and shows that the limit is described by the famous semi-circle law:
\begin{center}\includegraphics[width=4cm]{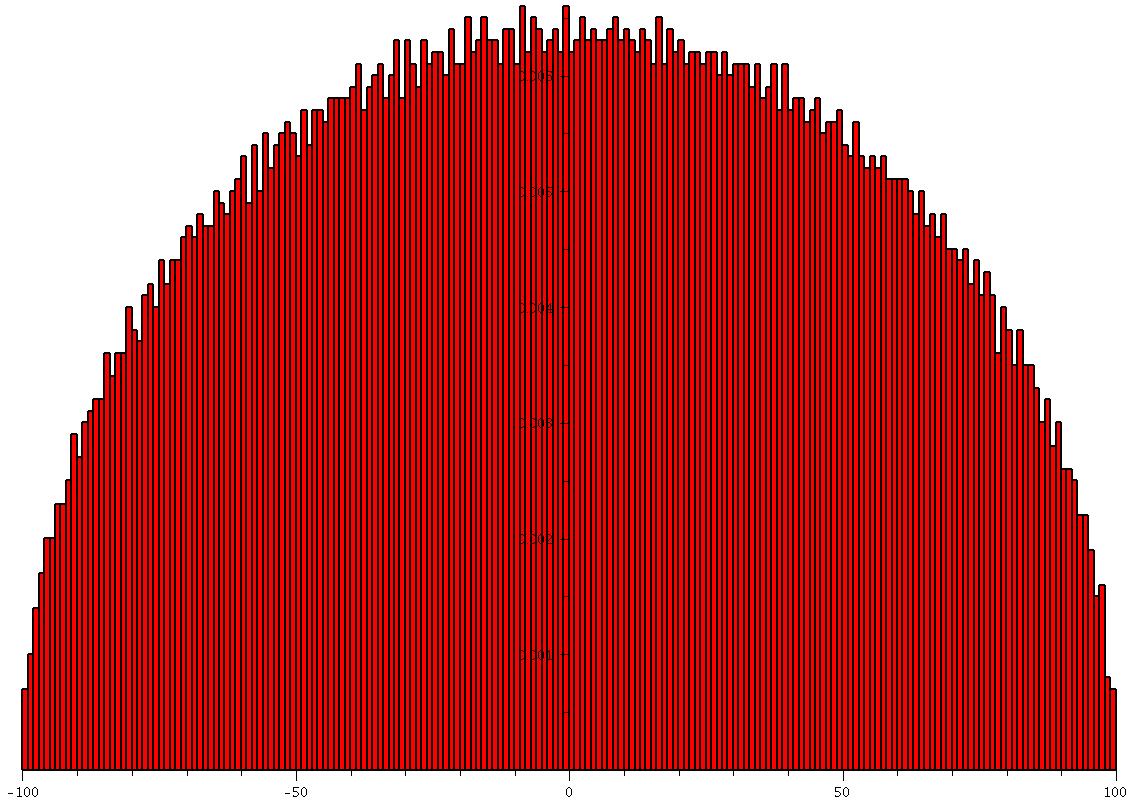}\end{center}

\begin{thm} [Wigner 56']\label{theowig}
Assume $pn$ goes to infinity as $n$ goes to infinity. Then, almost surely, for any $a<b$
$$\lim_{n\rightarrow\infty}\frac{1}{n}\#\{ i:\lambda_i^B\in \sqrt{np(1-p)}[a,b]\}=\lim_{n\rightarrow\infty}\frac{1}{n}\#\{ i:\lambda_i^X\in [a,b]\}=\sigma([a,b])$$
where $\sigma$ is the semi-circle law given by 
\begin{equation}
\sigma(dx)=\frac{1}{2\pi}\sqrt{4-x^2} 1_{|x|\le 2}
dx.
\end{equation}

\end{thm}
The semi-circle law is ubiquitous to random matrix theory as it describes the asymptotic behavior of random matrices with Gaussian entries but in fact any random matrix with independent centered entries $(a_{ij})_{i,j}$ such that
 $\E[|\sqrt{n}a_{ij}|^{2+\epsilon}]$ is uniformly bounded for some $\epsilon>0$. Such a convergence was proven first by Wigner in the case where $p$ is independent of $n$ based on the computation of the moments $\E[\tr \bX_n^k]$. Indeed, one can expand  the trace of moments of matrices in terms of the entries, and observe  that the indices which  contribute to the first order of this expansion can be described by  rooted trees, whereas $\sigma(x^k)$ is equal to the Catalan numbers which enumerate them. 

\subsection{Sparse Case}
On the other hand, the limiting distribution of the spectrum is very different when $pn$ is of order one. Namely we have the following theorem, see  \cite{KSV, ZAK}.
\begin{thm}\label{sparse}Assume $pn$ goes to $c\in (0,+\infty)$  as $n$ goes to infinity. Then, almost surely, for any $a<b$
$$\lim_{n\rightarrow\infty}\frac{1}{n}\#\{ i:\lambda_i^B\in \sqrt{np(1-p)}[a,b]\}=\lim_{n\rightarrow\infty}\frac{1}{n}\#\{ i:\lambda_i^X\in [a,b]\}=\mu_c([a,b])$$

\end{thm}
$\mu_c$ depends on $c$ and here are some plots :
\begin{center}\includegraphics[width=3.5cm]{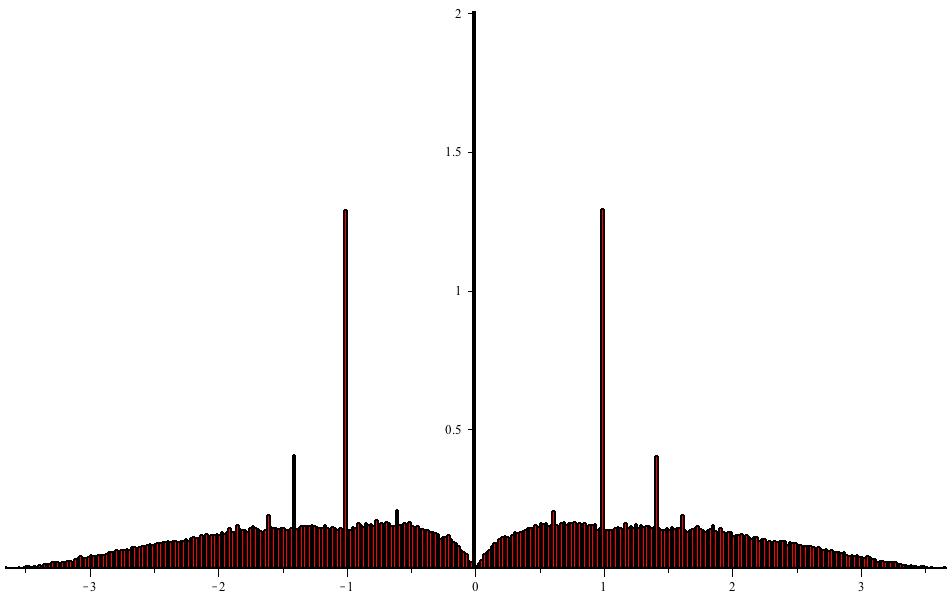}\hspace{0.5cm}\includegraphics[width=3.5cm]{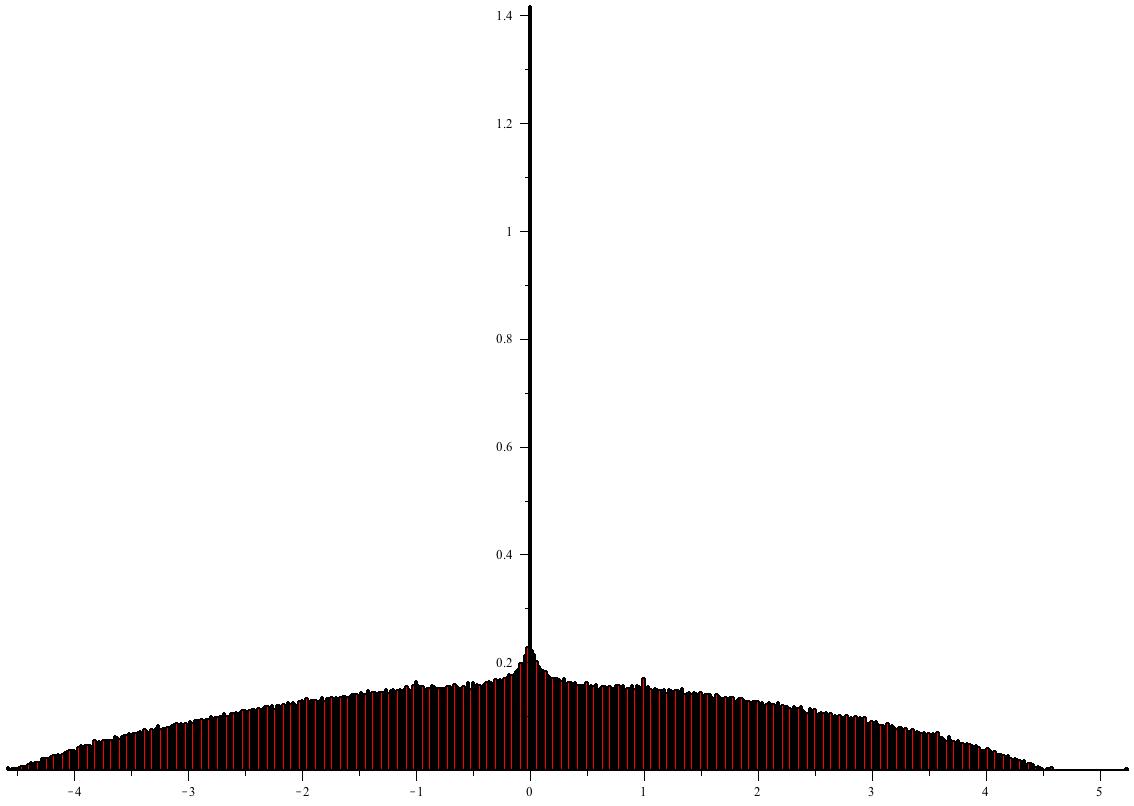}\hspace{0.5cm}\includegraphics[width=3.5cm]{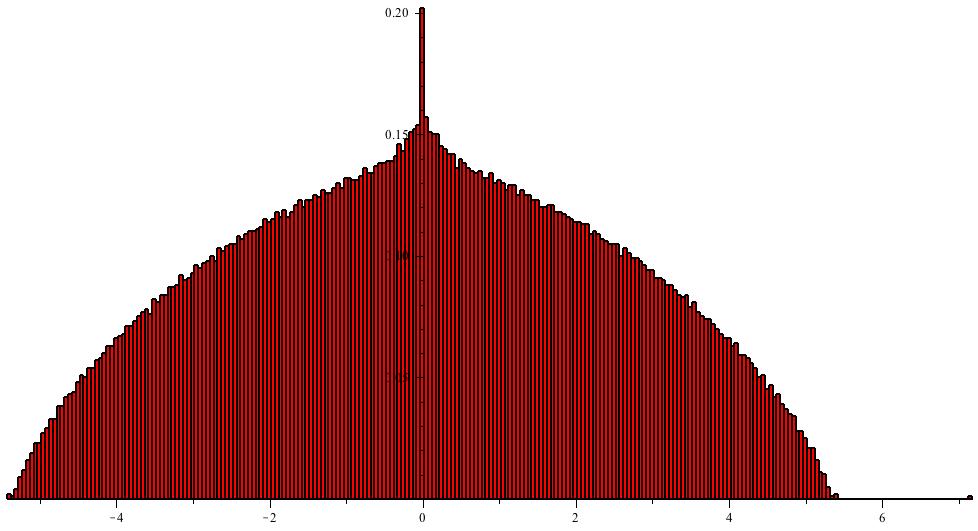}\end{center}

\centerline{Simulation for $c=1,2,3$ (Courtesy of J. Salez)}

The simulations indicate the presence of atoms. They were shown to be exactly given by totally real algebraic integers in \cite{Salez} for all $c>0$: these are the roots of monic polynomials with integer coefficients. It is easy to understand that the atoms should be totally algebraic integers as finite connected components are diagonal blocks with $0$ or $1$ entries whose characteristic polynomials have such roots. It is a much stronger statement to show that all such roots are atoms, in particular since totally algebraic integers are dense in the real numbers. $\mu_c$ has also 
a continuous spectrum: it was indeed proven in \cite{BSV} that $\mu_c$ has a non trivial continuous part if and only if $c>1$.This result is in fact hard to prove as the limit laws $\mu_c$ are described as the solution of complicated equations \cite{BoLe}, see also \cite{BAG1,BGM}.  However, such description could be used in \cite{ArBo} to prove the existence of an absolutely continuous part for sufficiently large $c$. Moreover, the first order expansion of $\mu_{c}$ in $c$ going to infinity was derived in \cite{EnMe}. The spectrum at the origin seems to have a Dirac mass whose weight could be computed \cite{BoLeSa}. 

\subsection{Idea of the proof}
The first proof of Theorem \ref{theowig} estimated the moments $\frac{1}{n}\tr(\bX_n)^k$ for all integer numbers $k$,  see \cite{Wig58} for the the first theorem and \cite{KSV, ZAK, BG08} for the sparse case. However, in order to go into more local results like the behavior of the eigenvectors or the local fluctuations, and as well to have more explicit formulas for the limit law, it is more convenient to study the resolvent. This path  can be used to study the asymptotics of the spectral measure of any self-adjoint matrix $\bX_n$ with independent entries modulo the symmetry constraint, and was generalized to study heavy tails matrices in \cite{KSV,BG08,BGM} based on the ideas  from \cite{BouchaudCizeau}. 
The idea is to derive the asymptotics  of the Stieljes transform 
$$G_n(z)=\frac{1}{n}\tr( z-\bX_n)^{-1}=\frac{1}{n}\sum_{i=1}^n\frac{1}{z-\lambda_i^X} $$ for a complex number $z$ away from the real line.

To this end, we use the Schur complement formula which reads
\begin{equation}\label{schur}(z-\bX_n)^{-1}_{ii}=\frac{1}{z-X_{ii}-\langle X_i, (z-\bX^{(i)})^{-1} X_i\rangle}\end{equation}
where $X_i=(X_{ij})_{j\neq i}$ and $\bX^{(i)}$ is the associated pinciple minor, namely the $(N-1)\times(N-1)$ matrix obtained from $\bX_n$ by removing the $i$th row and column. $X_{ii}$ goes to zero with $N$ and we can check (e.g by estimating the $L^{2}$ norm of the difference) that with probability going to one
\begin{equation}\label{app}
\langle X_i, (z-\bX^{(i)})^{-1} X_i\rangle= \sum_{j: j\neq i} X_{ij}^2(z-\bX^{(i)})^{-1}_{jj}+o(1)\,.\end{equation}

This is were the `light tail'' hypothesis $pn$ going to infinity starts to matter. Then, the entries $X_{ij}^2$ go to zero and have variance $1/n$ so that, since the $X_{ij}$ are independent of $\bX^{(i)}$, the law of large numbers (or a second moment computation) asserts that with probability going to one 
$$\sum_{j: j\neq i} X_{ij}^2(z-\bX^{(i)})^{-1}_{jj}=\sum_{j: j\neq i} \E[X_{ij}^2](z-\bX^{(i)})^{-1}_{jj}]+o(1)= \frac{1}{n} \sum_{j: j\neq i} (z-\bX^{(i)})^{-1}_{jj}+o(1)\,.$$
But again $\bX^{(i)}$ and $\bX_n$ vary only by a rank two matrix (if we complete $\bX^{{(i)}}$ by zero entries at the $i$th row and column), so that their spectrum is interlaced by Weyl's interlacing property. As a consequence
$$\frac{1}{n} \sum_{j\neq i} (z-\bX^{(i)})^{-1}_{jj}=\frac{1}{n} \sum_{i} (z-\bX_n)^{-1}_{jj}+O(\frac{1}{\Im (z) n})\,.$$
This approximation, together with \eqref{schur} and \eqref{app} implies that with high probability
\begin{equation}\label{quad}G_n(z)=\frac{1}{n} \sum_{i} (z-\bX_n)^{-1}_{jj}=\frac{1}{z-G_n(z)}+o(1)\,.\end{equation}
After recalling that $G_n(z)$ goes to zero as $N$ goes to infinity, we conclude that  since $G_n(z)$ goes to zero as the imaginary part of $z$ goes to infinity:
$$G_n(z)=\frac{1}{2}(z-\sqrt{z^2-4})+o(1)$$
is approximately the Stieljes transform of the semicircle law $G_\sigma(z)=\frac{1}{2}(z-\sqrt{z^2-4})$. Since $G_n$ is analytic and uniformly bounded for $\Im z>\epsilon$, Montel's theorem implies that $G_n$ converges to this limit away from the real line, which yields the vague convergence of the empirical measure of the eigenvalues. Because  $\frac{1}{n}\tr(\bX_n^2)$ is in $L^1$, the weak convergence follows.

On the contrary, in the heavy tails case  where $pn$ is of order one, the entries of $X_{ij}$ are often very small but of order one with a positive probability. Hence, the previous law of large numbers does not hold true any more and we can not expect such a simple equation as \eqref{quad}. 
In fact
$\sum_{j\neq i} X_{ij}^2(z-\bX^{(i)})^{-1}_{jj}$, if it converges, will a priori converge to a random variable. To study this convergence, we make the following assumption on the law $\mu_n$ of $X_{ij}$: 
\begin{equation}\label{Fou}\lim_{n\ra\infty}n\left(\mathbb \int (e^{{-}iu x^2}-1)d\mu_n(x)\right)=\Phi(u)\end{equation}
with $\Phi$ such that there exists $g$ on $\mathbb R^+$, with $g(y)$ bounded by $C y^\kappa$ for some $\kappa>-1$, such that for $u\in\bC^-$, 
\begin{equation}\label{hypPhi}\Phi(u)=\int_0^\infty g(y) e^{\frac{iy}{u}}dy. \end{equation}
This is satisfied  by the adjacency matrix of Erd\"os-R\'enyi graph with $\Phi(u)=c(e^{iu}-1)$ if $pn$ goes to $c$ and $g$  is a Bessel function \cite{BGM}, but also for other cases, for instance
for  $\alpha$ stable laws with $\Phi(u)=c(iu)^{\alpha/2}$ and $g(y)=C y^{\alpha/2-1}$ for some constants $c,C$. 
Then, it was shown in \cite{BG08, BGM} that $G_n(z){=\frac{1}{n}\tr(z-\bX^n)^{-1}}$ converges almost surely towards
$G$  given by\begin{equation}\label{St}G(z)=i \int e^{itz} e^{\rho_z(t)} dt, \,  z\in \mathbb C^+
\end{equation}
where $\rho_z:\bR^+\ra\{x+iy; x\le 0\}$ is the unique solution, analytic in 
$z\in\bC^+$, of the non linear equation
\begin{equation}\label{eqrho}\rho_z(t)=\int_0^\infty g(y) e^{\frac{iy}{t}z +\rho_z(\frac{y}{t})} dy\,.\end{equation}
This entails the convergence of the spectral measure of $\bX^n$, with $\sigma$ replaced by a probability measure with   Stieltjes transform given by \eqref{St}. The argument  to prove \eqref{St} and \eqref{eqrho} is as follows. We first remark that $G_n$ concentrates in the sense that it is close to its average, see Theorem \ref{concber}. We let $\rho^n$ be the  order parameter $\rho_z^n(x):=\mathbb E[\frac{1}{n}\sum  \Phi(x(z-\bX^{(i)})^{-1}_{jj})]$.
By \eqref{schur} and \eqref{app}, we find that if $\Im z>0$
\begin{eqnarray*}
G_n(z)&\simeq& \E[G_n(z)]=-i \E\left[\int_0^\infty e^{it z- it \sum_{j\neq i} X_{ij}^2(z-\bX^{(i)})^{-1}_{jj}} dt \right] +o(1)\\
&=& i \int_0^\infty e^{it z}\E\left[\prod_{j\neq i} 
\E[ e^{-it  X_{ij}^2(z-\bX^{(i)})^{-1}_{jj}}] dt \right]+o(1) \\
&=&-i \int_0^\infty e^{it z}\E\left[\prod_{j\neq i} \left(1+\frac{1}{n} \Phi( t(z-\bX^{(i)})^{-1}_{jj}))
\right) dt \right]+o(1) \\
&=&i \int_0^\infty e^{it z+\rho^n_z(t)} dt +o(1)\end{eqnarray*}
To conclude, we  need to show the convergence of $\rho^n$.  But $\rho^n$ can be seen to be analytic away from the real axis, and uniformly bounded under our hypothesis. This is enough to see that it is tight and any limit point will be analytic  by Montel theorem. Hence, it is enough to show that it has a unique limit point for $z$ with large imaginary part. To this end, we 
get an equation for $\rho^n$ which follows from \eqref{hypPhi} by
\begin{eqnarray*}
\rho^n_z(t)&=& \int_0^\infty g(y) \E[ e^{\frac{iy}{ x(z-\bX^{(i)})^{-1}_{11}}}]dy\\
&\simeq&  \int_0^\infty g(y) \E[ e^{\frac{iy}{ x}(z- \sum_{j\ge 2} X_{ij}^2(z-\bX^{(1)})^{-1}_{jj})}] dy+o(1)\\
&\simeq&  \int_0^\infty g(y) e^{\frac{iy}{ x}z}e^{ \rho^n_z(\frac{y}{x})} dy +o(1)\end{eqnarray*}
where in the second line we used \eqref{schur} and \eqref{app}.   One can conclude by proving  the uniqueness of the solutions to this equation when $z$ is far from the real line by showing that the non linear equation is then a contraction. 
The above arguments were made complete in \cite{BG08, BG14,BGM}. Another approach to heavy tails matrices and sparse Bernoulli matrices based on the Aldous 's PWIT was proposed in \cite{BCC}. 

\subsection{Extreme eigenvalues}

The asymptotic behavior of the extreme eigenvalues also depend on $c$ : they stick to the bulk when $pn\gg \ln n$ and then go away at distance of order $\sqrt{ \ln n}$. We more precisely have the following result, putting together the article of 
Benaych-Georges, Bordenave, Knowles \cite{BCK1} and Alt, Ducatez, Knowles \cite{ADK1}, see also \cite{TiYo}.
\begin{thm}\label{asmax}\begin{itemize}
\item Assume {{$pn/\ln n\ra+\infty$}}, then the largest eigenvalue of $\bX_{n}$ sticks to the bulk : {{$\lambda_{1 }^X\ra 2$ }}.
\item Assume {{$pn/\ln n\ra 0$}}, then  {{$\lambda_{1}^X\simeq \sqrt{\ln n/\ln (\ln n/pn)}$}}.
\item Assume {{$pn\simeq C\ln n$}}, then for ${{C>1/(\ln 4-1):=C^{*}}}$ the eigenvalues stick to the bulk, whereas for ${{C<1/(\ln 4-1})}$
$${{\lambda_{1}^X=\frac{\alpha}{\sqrt{\alpha-1}}, \alpha=\max\frac{1}{pn}\sum_{j} B_{ij}}}\,.$$
\end{itemize}

\end{thm}
Observe that $\sum_{j} B_{{ij}}$ is the degree of vertex $i$: the largest eigenvalue is hence created by the largest degree in the graph. 
In fact, in Alt, Ducatez, Knowles \cite{ADK1}, it is shown  that all eigenvalues outside the bulk are created by  vertices with large degrees 
 when $pn\le C^{*}\ln n$. 
\section{Fluctuations}
\subsection{Concentration of measure}
Concentration of measure has become a central tool in probability and in particular in random matrix theory. It allows us to prove that some quantities such as  smooth function of independent variables, are not much random. It was crucial in the previous proof of the convergence of the spectral measure. However, it generally depends on the tails of the random variables.  Herbst's argument allows considering random variables with sub-Gaussian tails and more precisely random variables whose distribution satisfies log-Sobolev inequalities, which is the case for instance when their density is strictly log-concave as for Gaussian's variables. To deal with bounded variables such as  the entries of Bernoulli matrices, one should rather use the theory developed by Talagrand \cite{Tala}. This was done in \cite{GZconc} where the spectrum of random matrices was observed to be a smooth function of its entries and the associated Lipschitz norm was computed. It resulted in the following theorem \cite[Theorem 1.1]{GZconc}. We hereafter consider a symmetric matrix $\bA$ with independent entries above the diagonal with distribution $a_{ij}/\sqrt{n}$ where $a_{ij}$ is distributed according to  $P_{ij}$  supported in a compact set $K$ with width $|K|$. 
\begin{thm} \begin{enumerate}
\item Take $f$  convex and Lipschitz with Lipschitz norm $\|f\|_L$. Then, for any $\delta>\delta_0(n)= 8|K|\|f\|_L/n$, 
$$\Pp\left(| \frac{1}{n}\tr(f(\bA))- \E[\frac{1}{n}\tr(f(\bA))]|>\delta\|f\|_L
\right)\le 4 \exp\{-n^2\frac{(\delta-\delta_0(n))^2}{16|K|^2}\} $$
\item There exists a finite constant $c>0$ such that for any $\delta> \delta_1(n)\simeq \sqrt{\delta_0(n)}$
$$\Pp\left(\sup_{f\in \mbox{Lip}_{\mathcal K} } | \frac{1}{n}\tr(f(\bA))- \E[\frac{1}{n}\tr(f(\bA))]|>\delta\|f\|_L\right)\le \exp\{-n^2\frac{(\delta-\delta_1(n))^2}{c|K|^2} \}\,.$$

\item Let $\lambda_1^A$ be the largest eigenvalue of $\bf A$. Then
$$\Pp\left( |\lambda_1^A-\E[\lambda_1^A]|\ge  \delta |K| \right)\le \exp\{-\frac{(\delta- 8|K|/\sqrt{n})^2n }{16 }\}\,.$$
\end{enumerate}
\end{thm}
This result is a direct application of Talagrand's beautiful theory and the computation of Lipschitz constants of functions of the spectral measure in terms of the entries, see \cite{GZ3,AGZ}. The original statement proves concentration around the median rather than the mean, but it is easy to go from one result to the other up to some error $\delta_{0}(n), \delta_{1}(n)$. The second point is deducted from the first by approximating a general function by convex functions.
It applies to Bernoulli matrices straightforwardly by taking $|K|=1/\sqrt{p(1-p)}$
\begin{thm}\label{concber} Take $f$  convex and Lipschitz with Lipschitz norm $\|f\|_L$. Then, for any $\delta>\delta_0(n)= 8\sqrt{\pi}|f|_L /n p(1-p)$, 
$$\Pp\left(| \frac{1}{n}\tr(f(\bX_n))- \E[\frac{1}{n}\tr(f(\bX_n))]|>\delta+\delta_0(n)\right)\le \exp\{-p(1-p)n^2\frac{(\delta)^2}{16|f|_L^2}\}\,. $$
Moreover, for any $\delta>\delta_0'(n)=O(1/\sqrt{p(1-p)n})$
$$\Pp\left(| \lambda_1- \E[\lambda_1]|>\delta+\delta_0'(n) \right)\le \exp\{-p(1-p)n\delta^2\} $$
\end{thm}
As we can see, the speed of the concentration deteriorates with $p$ going to zero to be of order $n$ when $np$ is of order one. In fact, it can be shown that the worse concentration estimates for the empirical measure are of the order of exponential in $n$. In fact, we have the following result due  to C. Bordenave, P. Caputo and D. Chafai
\cite{BCC} which is
based on the Azuma's-Hoeffding inequality and requires only the independence of the vectors of the random matrix.
\begin{lem}\label{BCClemma}
Let $\|f\|_{TV}$ be the total variation norm,
$$\|f\|_{TV}=\sup_{x_1<\cdots<x_p}\sum_{i=2}^p|f(x_i)-f(x_{i-1})|
$$
Then, for any self-adjoint matrix $\bX_n$ with independent vectors $((X_{ij}, i\le j),1\le j\le n)$ and eigenvalues $(\lambda_i)_{1\le i\le n}$
for  any function $f$ with finite total variation and any $\delta>0$,
norm so that $E[|\frac{1}{n}\sum_{i=1}^n f(\lambda_i)|]<\infty$,
$$P\left(\left|\frac{1}{n}\sum_{i=1}^n f(\lambda_i)-\E[\frac{1}{n}\sum_{i=1}^n f(\lambda_i)]\right|\ge \delta\|f\|_{TV} \right) \le 2
e^{-\frac{n\delta^2}{8 }}$$
\end{lem}
In the general case however, the extreme eigenvalues do not concentrate and can be very large for heavy tails entries \cite{ABP, ADK1}.
\subsection{Global fluctuations}
It is a natural question to wonder how the empirical measure of the eigenvalues fluctuates and in particular whether the concentration result of Theorem \ref{concber} is on the optimal scale. 
In the case where $p$ is of order one, this question was first answered by Jonsson \cite{jonsson} by estimating moments, and 
in the context of Gaussian matrices by Johansson \cite{johansson} by using loop equations. The main point is that the central limit theorem does not require a renormalization by the famous $\sqrt{n}$ as for the classical central limit theorems.
\begin{thm} Assume $p\in (0,1)$ independent of $n$. Let $f$ be a continuously differentiable function. Let $\lambda_i$ be the eigenvalues of $\bX_n$. Then
$$\sum_{i=1}^n f(\lambda_i)- \E[\sum_{i=1}^n f(\lambda_i)]$$ converges in distribution towards a centered Gaussian variable with variance
$$V(f)=\frac{1}{2\pi^2}\int_{-2}^2\int_{-2}^2\left(\frac{f(x)-f(y)}{x-y}\right)^2 \frac{(4-xy)}{\sqrt{4-x^2}\sqrt{4-y^2}}dx dy\,.$$ \end{thm}
The central limit theorem also holds if one recenters with respect to the limit rather than the expectation, 
 see e.g. \cite{PaSh}.

On the contrary, if $pn$ goes to a constant $c$ we see that Theorem \ref{BCClemma} gives the optimal speed and we have  a ``more'' classical central limit theorem \cite{ShTi,BGM, AZ08}:
\begin{thm}\label{cltsparse} Assume that $pn$ goes to $c\in (0,+\infty)$. Let $f$ be a $C^1_b$ function.  Then
$$\frac{1}{\sqrt n}\left( \sum_{i=1}^n f(\lambda_i)- \E[\sum_{i=1}^n f(\lambda_i)]\right)$$ converges in law towards a centered Gaussian variable with non trivial variance.

\end{thm}

Together with \cite{HeYu}, we claim that at least for $pn$ of order one, or in $[n^{\varepsilon},n^{1-\varepsilon}], $ or $p$ of order one, we have (probabilibly for all $p>1/n$), 
\begin{thm}\label{cltsparse} Let $f$ be a $C^1_b$ function.  Then
$${\sqrt p}\left( \sum_{i=1}^n f(\lambda_i)- \E[\sum_{i=1}^n f(\lambda_i)]\right)$$ converges in law towards a centered Gaussian variable with non trivial covariance.
\end{thm}
This result should hold for any $p>1/n$.

\subsection{Local laws}
An important breakthrough toward the understanding of  local fluctuations and eigenvectors is to analyze the so-called local laws as foreseen in \cite{ESY10}. Namely to estimate $\sum f(\lambda_i)$ for less smooth functions, in fact for functions on a mesoscopic scale $f(x)=g(N^\alpha(x-E))$ for some $\alpha\in (0,1)$. Equivalently, one can look at $f(x)=(z-x)^{-1}$ with $z=E+i \eta$ with $\eta$ of order $N^{-\alpha}$ (indeed the later can serve to approximate conveniently the first). In this scale, it was proven that if $pn$ goes to infinity, the mesoscopic distribution of the eigenvalues is still very close from the semi-circle distribution. 
Indeed, let us define  the Stieltjes transform to be given by 
$$G_n(z)=\frac{1}{n}\sum_{i=1}^n\frac{1}{z-\lambda_i}, G_\mu(z)=\int \frac{1}{z-\lambda}d\mu(\lambda)$$
 In \cite[Theorems 2.8 and 2.10]{EKYY}, the following result was proven, where $\zeta$- high probability meaning  a probability greater or equal to $1-e^{-v(\ln n)^\zeta}$ with for some $v>0$. 
\begin{thm}\label{theohy} 
There are universal constants $C_1,C_2>0$ such that the following holds. Assume that
$$pn\ge (\log n)^{C_1\xi}, \xi= C_2\log \log n$$
Then, for $E\in [-3,3]$ and  $D=\{z=E+i\eta, 0<\eta<3\}$, 
$$\cap_{z\in D}\left\lbrace |G_n(z)-m_\sigma(z)|\le (\log n)^{C_2 \xi}\left(\min\{\frac{1}{pn\sqrt{\kappa_E+\eta}},\frac{1}{\sqrt{pn}}\}+\frac{1}{n\eta}\right)\right\rbrace$$
holds with $\zeta$-high probability. Moreover, for $\eta>(\ln n)^{C\zeta} n^{-1}$
$$\#\{i:\lambda_i\in [E-\eta,E+\eta]\}=n\sigma([E-\eta,E+\eta])\left(1+O(\ln n)^{C\zeta}\left(\frac{1}{n\eta^{\frac{3}{2}}}+\frac{1}{pn\eta}\right)\right)$$
with  $\zeta$-high probability.
\end{thm}
The above theorem applies for any $p$ such that $pn$ goes to infinity much faster than any $\ln n$, see e.g. \cite{ADK1}.  Below $\ln n$, the extreme eigenvalues were shown to be dictated by the largest degree in the graph \cite{BCK1}.

A similar statement in the sparse case where $pn$ goes to a finite constant is still open. Indeed, the fact that $\mu_c$ has a dense set of atoms and a continuous part makes the analysis a priori much more involved and the local law more difficult to conjecture.
An easier heavy tail matrix model was studied in \cite{BG08,BouchaudCizeau,BG}, namely the random matrices with alpha-stable independent entries. In this case, the entries follow the alpha-stable law $\mathbb P\left(|A_{ij}|\ge t\right)\simeq t^{-\alpha}/n$. When $\alpha<2$, it was shown in \cite{BG08,BouchaudCizeau} that the empirical measure converges towards a limiting law $\mu_\alpha$ which is different from the semi-circle law.  One of the advantages of this model is that $\mu_\alpha$ is absolutely continuous except possible for a discrete set of atoms. 
Of course, one can not expect the eigenvalues to be as rigid in the heavy tails case since this would contradict the central limit theorem (which holds as in Theorem \ref{cltsparse}, see \cite{BGM}). Hence, in this case, large eigenvalues should be less rigid, creating large fluctuations. The following result was proved if the $A_{ij}$ are $\alpha$-stable variables in \cite{BG,BG2}:
 for all $t \in \bR$, 
\begin{equation}
\E [\exp ( i t A_{11} ) ] = \exp ( - \frac{1}{n} w_\alpha |t |^\alpha ),
\end{equation}
for some $0 < \alpha < 2$ and $w_\alpha = \pi /  ( \sin(\pi \alpha / 2)\Gamma(\alpha)) $. We  put 
\begin{equation}
\rho = \left\{ \begin{array}{lcl}
\frac 1 2 & \hbox{ if } & \frac 8 5 \leq \alpha < 2 \\
\frac \alpha { 8 - 3 \alpha } & \hbox{ if } & 1 < \alpha < \frac 8 5 \\
\frac {\alpha} { 2+ 3 \alpha } & \hbox{ if } & 0 < \alpha \leq 1.
\end{array} \right.
\end{equation}
Then, there exists a finite set ${\mathcal E}_\alpha \subset \bR$ such that if $K \subset \bR \backslash {\mathcal E}_\alpha $ is a compact set and $\delta >0$, the following holds. There are constants $c_0, c_1 >0$ such that for all integers $n \geq 1$, if $I\subset K$ is an interval of length $ |I|\ge c_1 n^{-\rho} (\log n)^2$, then
\begin{equation}
 \left|  N_I -n \mu_\alpha ( I) \right| \leq \delta n  |I |,
\end{equation}
with probability at least $1 - 2 \exp \left( - c_0 n \delta^2 | I | ^2 \right)$. The fact that our result might not be true on a finite set of values should only be technical . 
This result was improved in \cite[Theorems 3.4 and 3.5]{ALY} in order to tackle $I$ of size o $n^{-\omega(\alpha)}$ with $\omega(\alpha)>1/2$ (and $\Re(z)$ small enough when $\alpha<1$). Such an optimal scale is important in the of study the local fluctuations of the spectrum. 

In both light and heavy tails, the main point is to estimate the Stieltjes transform $G_n(z)=\frac{1}{n}\sum_{i=1}^n (z-\lambda_i)^{-1}
$ for $z$ going to the real axis : $z=E+i\eta$ with $\eta$ of order  nearly as good as $n^{-1}$ for light tails, $n^{-\rho}$ for  heavy tails. This is done by showing that 
 $G_n$ is characterized approximately by a closed set of equations.
In the case of lights tails, one has simply a quadratic equation for $G_n$ and needs to show that the error terms remain small as $z$ approaches the real line. In the heavy tails case, the equations are much more complicated, see \eqref{St} and \eqref{eqrho}, and therefore more difficult to handle. Similar questions are completely open for other heavy tails matrices, including Bernoulli matrices with $pn$ of order one.

\subsection{Local fluctuations}
When the average degree $pn$ is large, one expects the eigenvalues to behave exactly as the eigenvalues of a symmetric matrix with independent Gaussian entries. The so-called GOE matrices. The advantage of Gaussian matrices is that they are an integrable model of random matrices in the sense that many of their properties  can be exactly computed. To start with, the joint distribution of its eigenvalues $(\lambda_{i}^{G})_{1\le i\le n}$ is explicit:
\begin{equation}\label{lawG}d\Pp(\lambda^{G})=\frac{1}{Z} \Delta(\lambda) e^{-\frac{n}{4}\sum( \lambda_i^{G})^2}\prod d\lambda_i^{G}\end{equation} 
where $\Delta(\lambda)=\prod_{i<j} |\lambda_i-\lambda_j|$ is the Vandermonde determinant. In particular, this formula does not depend on the eigenvectors. Based on this formula, Tracy and Widom could study the local fluctuations 
of the spectrum  $(\lambda_i^G)_{1\le i\le n}$ \cite{TW1,TW5} and they  proved that
$$\lim_{n\ra\infty}\Pp( n^{2/3}(\lambda_{1}^G-2)\le s)=F_1(s)$$
where $F_1$ is the distribution function of the Tracy-Widom law. For the eigenvalues in the bulk, it was proven \cite{ME} that for all smooth compactly supported function
$$
\mathcal E_{\bG_{n}}(O,E)= \E[ O(n(\lambda_i^G-E),\cdots, n(\lambda_{i+p}^G-E)) ]$$
converges as $n$ goes to infinity and the limit is described in terms of Pfaffian distributions.

The universality in the bulk was obtained after the a series of works including notably \cite{TV, ESY}  in \cite[Theorems 2.5]{EKYY11} (for $\phi\ge 2/3$) and  improved in \cite{HLY} (for $\phi>0$) to finally get: 
\begin{thm}\label{bu}
(Bulk universality) Suppose $pn> n^{\phi}$ with $\phi>0$, there exists $b_n$ going to zero so that for all smooth compactly supported function $O$,  any $E\in (-2,2)$

$$\lim_{n\ra\infty} \int_{E-b_n}^{E+b_n} \frac{dE'}{2b_n}\left( \mathcal E_{\bG_{n}}(O,E')-\mathcal E_{\bB_{n}}(O,E')\right)=0$$
\end{thm}

Moreover, the universality at the edge was obtained in \cite[Theorems 2.7]{EKYY11}, see also \cite{soshnikov2},
\begin{thm}(Edge universality)
Suppose that $pn>n^{\phi}$, $\phi>2/3$. 
Then there exists $\delta>0$ such that
$$\Pp( n^{2/3}(\lambda_{2}^B-2)\le s)=\Pp( n^{2/3}(\lambda_{2}^G-2)\le s+O(n^{-\delta}))+O(n^{-\delta})$$
\end{thm} 
This statement was generalized to $pn>n^{1/3}$ but the largest eigenvalue then needs to be shifted by a deterministic drift of order $1/pn$ \cite{LS}. 
Beyond this threshold, the fluctuations of the second largest eigenvalue starts to be Gaussian.

When $pn$ decreases below $1/3$, it was proven that universality stops to hold and fluctuations of the largest eigenvalue start to be Gaussian. The precise transition between Tracy-Widom law and Gaussian fluctuations when $p$ is of order $n^{-2/3}$ was described  \cite{HLY}. When $n^{o(1)}\ll pn \ll n^{1/3}$, \cite{HLY, KY} shows that the fluctuations of the extreme eigenvalues are Gaussian, even if they stick to the bulk. In the case where $pn\ll \ln n$,
Theorem \ref{asmax} asserts that the eigenvalues go away from the bulk, at distance of order $\sqrt{\ln n}$. 
The corresponding eigenvectors are localized close to the vertices with a high degree. In an even more recent preprint \cite{ADK3}, the same authors show that these eigenvalues follow a Poisson point process. 

Such questions are open for Bernoulli random matrices  with $pn$ of order $c\in (0,+\infty)$ and eigenvalues in the bulk. Indeed, as we have seen, the limiting density is a mixture of atoms and continuous density and it is not yet clear how to zoom in the spectrum in such a situation. 
However, such questions could be analyzed for L\'evy matrices with  $\alpha$-stable entries in the regime where local law can be obtained on the optimal scale $n^{-1/2}$ \cite{ALY}. In fact,
one expects the following transition to occur, see   \cite{TBT}:

$\bullet$ If $\alpha\in [1,2]$,
all eigenvectors are corresponding to finite eigenvalues are completely delocalized. Further, for any $E\in\mathbb R$ the local statistics of the eigenvalues near $E$ converge to those of the GOE as $N$ goes to infinity.

$\bullet$ If $\alpha\in (0,1)$, There exists a mobility edge $E_\alpha$ such that for $|E|<E_\alpha$    the local statistics of the eigenvalues near $E$ converge to those of the GOE as $N$ goes to infinity. But if  $|E|>E_\alpha$    the local statistics of the eigenvalues near $E$ converge to those of a Poisson point process and all eigenvectors in this region are localized.

\begin{center}
\includegraphics[width=7cm]{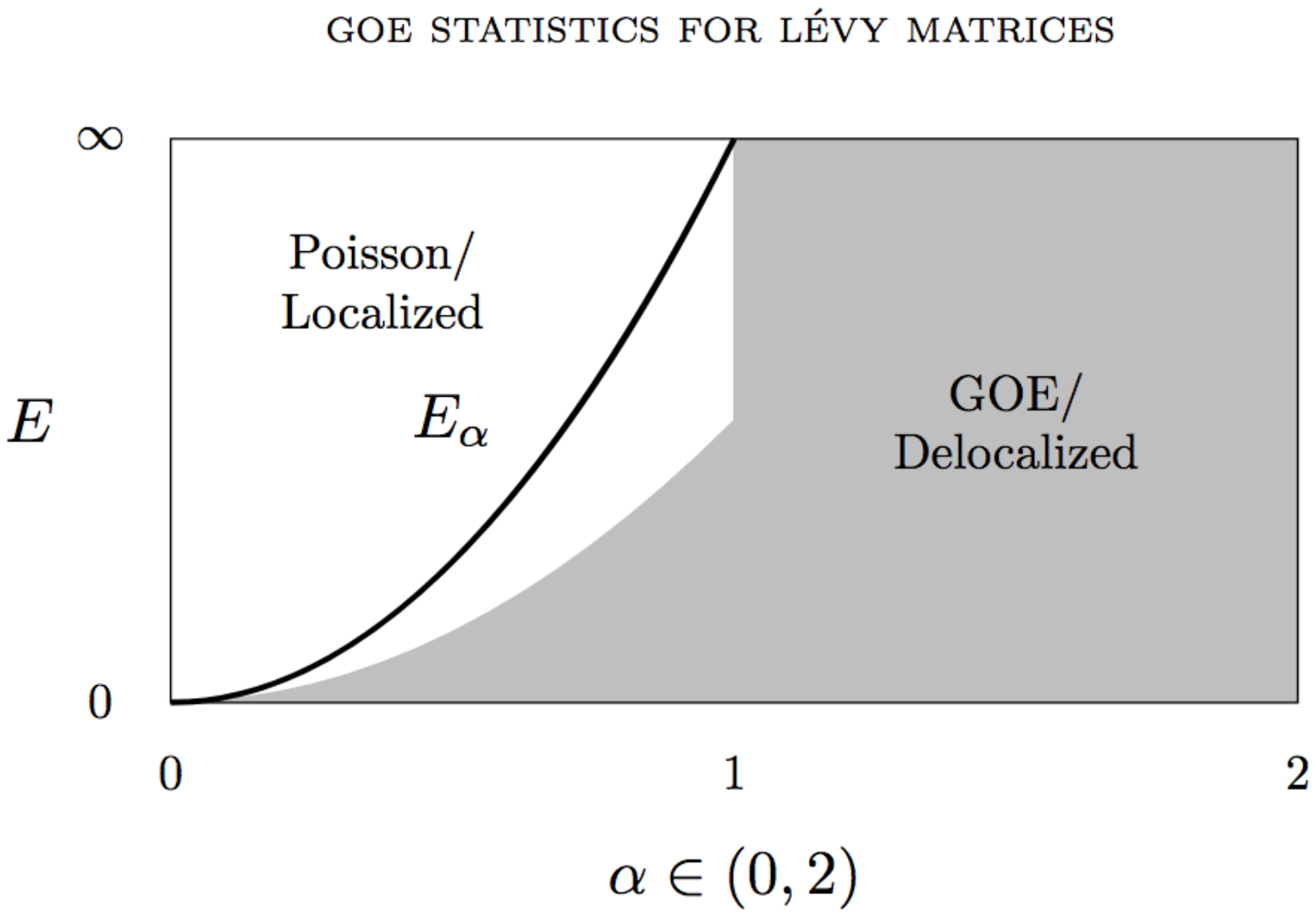}\end{center}

The fact that local statistics are given by those of Gaussian matrices for $\alpha\in (1,2)$ or $\alpha\in (0,1)$ and $E$ small enough, except for $E$ in some finite set was proven in \cite[Theorems 2.4 and 2.5]{ALY}.

\subsection{Properties of the eigenvectors}
The properties of the eigenvectors are intimately related with local laws. Indeed, by definition of the eigenvectors, if $v$ is an eigenvector of  the symmetric matrix $\bX_n$ for the eigenvalue $E$ and we set $\langle v, e_i\rangle=v_i$, $X_1$ the 
first column vector of $\bX_n$ while $\bX_n^{(1)} $ is the $(n-1)\times (n-1)$ principal minor of
$\bX_n$ obtained by removing the column and row vector given by $X_1$ and $X_1^T$

$$v_1^2=(1+\langle X_1, (E-\bX_n^{(1)})^{-2} X_1\rangle)^{-1}$$
where, at least in the dense cases $\langle X_1, (E-\bX_n^{(1)})^{-2} X_1\rangle$ is close to $\frac{1}{n}\tr(E-\bX_n)^{-2} $, and so is governed by the local law.
In \cite[Theorems 2.16]{EKYY}, it was proven that
\begin{thm}(Complete delocalization of eigenvectors)Assume the hypotheses of Theorem \ref{theohy} with  $pn>n^\phi$ with $\phi>0$. Let $v_i$ be the eigenvectors of $\bB_n$ for the eigenvalues $\lambda_n\le \lambda_{n-1}\cdots\le \lambda_1$. Then
$$\max_{i\le n}\|v_i\|_\infty\le \frac{(\ln n)^{4\zeta}}{\sqrt{n}}$$with $\zeta$-high probability. 
\end{thm}
This result was extended to $q$ going to infinity logarithmically only more recently \cite{ADK2}. We roughly state their result: 

\begin{itemize}
\item (Semilocalized phased) Assume $C\sqrt{\ln n}\ln\ln n\le \sqrt{pn} \le 3\ln n$ and $w$ be a normalized eigenvector of $\bB_n$ with non trivial eigenvalue $E\ge 2+C\zeta^{1/2}$. We let $\Lambda(\alpha)=\alpha/\sqrt{\alpha-1}$ and $\alpha_x=\sum_y \bB_{xy}/pn$. We let $W_{E,\delta}$ be the set of vertices such that $\Lambda(\alpha_x)\in[E-\delta,E+\delta]$. Then for each $x\in W_{E,\delta}$ there exists a normalized vector $v(x)$ supported in a ball around $x$ and radius $c\sqrt{\ln n}$, such that the support of $v(cx)$ and $v(y)$ is distinct if $x\neq y$ and 
$$\sum_{x\in W_{E,\delta}}\langle v(x),w\rangle^2\ge 1-C(\sqrt{\ln n}{pn}\ln pn +\sqrt{\ln n}{pn}\frac{1}{E-2})^2\delta^{-2}$$
Moreover $$\sum_{y\in B_r(x)}(v(x))^2_y\le \frac{1}{(\alpha_x-1)^{r+1}}$$
\item(Delocalized phase) For any $\nu>0$ and $\kappa>0$ there exists a constant $C>0$ such that for $pn\in [C\sqrt{\ln n}, (\ln n)^{3/2}]$, if $w$ is a normalized eigenvector for $\bB_n$ with eigenvalue $E\in [-2+\kappa,-\kappa]\cup[\kappa,2-\kappa]$,
$$\|w\|_\infty^2\le n^{-1+\kappa}$$
with probability greater than $1-n^{-\nu}$.
\end{itemize}
This question is completely open for Bernoulli random matrices with $pn$ of order one  but the understanding of L\'evy matrices is again more complete. Based on \cite{BG,BG2,ALY} we can assert that  Tarquini, Biroli, and Tarzia's conjecture   \cite{TBT} is partly proven. Indeed the complete delocalization is proven for $\alpha\in (1,2)$ and $\alpha\in (0,1)$ and small enough eigenvalues. A sort of localization for $\alpha\in (0,1)$ for large enough eigenvalue was derived in \cite{BG}, and was shown to be not true for small enough eigenvalues in \cite{BG2}: the transition and the value of the mobility edge is still an open question.  In fact, even in the case where the eigenvalue statistics belong to the universality class of Gaussian matrices, the fine properties of the eigenvectors of L\'evy matrices differ \cite{ALM}. Let us also mention \cite{RuVe} which shows under quite general assumptions that eigenvectors  are somehow uniformly delocalized in the sense that any subset of at least eight coordinates carries a non-negligible part of the mass of an eigenvector.

\section{Rare events}It is sometimes important to estimate the probability of rare events, such as the probability that the extreme eigenvalues take unlikely values or the empirical measure of the eigenvalues shows an unlikely profile, and what kind of optimal strategy can lead to such deviations from the expected behavior. In the case of Gaussian symmetric matrices, the joint density of the eigenvalues is known \eqref{lawG}.
 One finds by sort of Laplace's principle \cite{BAG, BDG} the large deviations for the empirical measure and the largest eigenvalue.
\begin{thm}
Let $\lambda_{n}^G\le\lambda_{2}^G\cdots\le \lambda_{1}^G$ be the eigenvalues of a GOE matrix. Then
\begin{itemize}
\item Let $E(\mu)=\frac{1}{2}\int\int (\frac{x^{2}}{4}+\frac{y^{2}}{4}-\ln|x-y| )d\mu(x)d\mu(y)$ and set $\mathcal E(\mu)=E-\inf E$. Then $\mathcal E$ is a good rate function and the distribution of the empirical measure of the eigenvalues $\hat\mu_{n}=\frac{1}{n}\sum \delta_{\lambda_{i}^G}$ satisfies a large deviation principle with speed $n^{2}$ with rate function $\mathcal I$, that is for every closed set $F$
$$\limsup_{n\rightarrow\infty}\frac{1}{n^{2}}\ln \mathbb P\left(\hat\mu_{n}\in F\right)\le-\inf_{F}\mathcal E$$
whereas for any open set $O$
$$\limsup_{n\rightarrow\infty}\frac{1}{n^{2}}\ln \mathbb P\left(\hat\mu_{n}\in O\right)\ge-\inf_{O}\mathcal E$$
\item Let $I_G(x)=\frac{1}{2}\int_{2}^{x }\sqrt{4-y^{2}} dy$ for $x\ge 2$ and $I_G(x)=+\infty$ for $x<2$. Then $I$ is a good rate function and the distribution of $\lambda_{1}^G$ satisfies a large deviation principle  with speed $n$ and  good rate function $I_{G}$.
\end{itemize}
\end{thm}
In this case, deviations of the spectrum can be created independently from the eigenvectors which stay uniformly distributed. 
On the other hand, if the entries have sharp exponential decay, large deviations can be created by large entries. 
Assume  that  for some $\alpha\in (0,2)$, there exists $a>0$ so that for all $i,j$
$${\lim_{t\rightarrow\infty } 2^{-1_{i=j}} t^{-\alpha}\ln \mathbb P( |\sqrt{n} X_{ij}|\ge t)=- a }$$
\begin{thm}
\begin{itemize}
\item  \cite{BordCap} The law of the empirical measure satisfy a LDP in the speed ${n^{1+\frac{\alpha}{2}}}$ and good rate function which is infinite unless $\mu=\sigma\boxplus \nu$ and then equals $a\int |x|^{\alpha}d\nu(x)$.
\item  \cite{fanny} The law of the largest eigenvalue satisfies a LDP with rate ${ n^{\frac{\alpha}{2}}}$ and GRF proportional to 
$(\int (x-y)^{-1}d\sigma(y))^{-\alpha}$.
\end{itemize}
\end{thm}

However, the situation is much less understood for Bernoulli matrices  and again the sparse and the dense regime lead to very different results and techniques.
We discuss these questions hereafter.

\subsection{Large deviations for the extreme eigenvalues}
 Let us first consider the dense case.  In \cite{HuGu,AGH}, we considered the large deviations for the largest eigenvalue of Wigner matrices and showed that if the entries are Rademacher, then the same large deviation principle holds, whereas in general there is a transition between deviations close to two where the rate function is  the Gaussian one whereas for large deviations towards large enough values the rate function is more of a heavy tail type. 
In a work in progress with F. Augeri, R. Ducatez and J. Husson, we prove that

\begin{thm}
\begin{itemize}
\item Assume that $p=1/2$. Then the law of $\lambda_1^X$ satisfies a large deviation principle in the scale $n$ and with the same rate function $I_G$ than for the GOE matrix.
\item Assume $p\in (0,1/2)$. Then for $x$ close enough to $2$, the probability that $\lambda_1^X$ is close to $x$ is the same than in the Gaussian case. But for $x$ large enough, 
$$\limsup_{\delta\downarrow 0}\limsup_{n\rightarrow \infty} \frac{1}{n}\ln \mathbb P(|\lambda_1^X-x|<\delta)=-I_p(x)$$
where $I_p(x)<I_G(x)$.
\end{itemize}
\end{thm}
The case $p\in (1/2,1)$ is under investigation. In fact, analyzing the large deviation requires to understand good strategies to create the deviations. For $p=1/2$ it is shown that an optimal strategy is to  tilt the law of the entries in order to  change their expectation so that the matrix looks like a rank one deformation of Bernoulli matrix with a delocalized deformation. The eigenvectors also stay delocalized through this deformation.  When $p<1/2$ and $x$ is large, it turns out that the optimal strategy is to create fully connected components of size $\sqrt{n}$. For $p>1/2$ the picture is less clear and we suspect that vertices with high degree are optimal ways to create large eigenvalues. 

Let us now consider the sparse case following 
 \cite{BBG}: in this case we already saw that large eigenvalues are created by  vertices with large degree, namely with row or column vectors with many entries equal to one.  
 \begin{thm} Let $L_{p}=\frac{\ln n}{\ln \ln n-\ln (np)}$ and assume $\ln(1/np)\ll\ln n$ and $np\ll \sqrt{\ln n/\ln\ln n}$. Let $\lambda_{2}$ be the second largest eigenvalue of $\bB_{n}$.
 Then for any $\delta\ge 0$, 
 $$\lim_{n\ra\infty}\frac{-\ln P\left(\lambda_{2}\ge (1+\delta) \sqrt{L_{p}}\right)}{\ln n}=2\delta+\delta^{2}$$
 whereas
 $$\lim_{n\ra\infty}\frac{-\ln P\left(\lambda_{2}\le (1-\delta) \sqrt{L_{p}}\right)}{\ln n}=2\delta-\delta^{2}$$
\end{thm}
 
\subsection{Large deviations for the empirical measure}
In \cite[Theorem 1.6]{BoCa15}, a large deviation for the empirical measure of the eigenvalue in the sparse case was derived : we do not precise the rate function as it is obtained by contraction from the large deviation for the empirical neighborhood distribution.
\begin{thm}Assume $pn$ is fixed. Then the law of $\hat\mu_{n}$ satisfies a large deviation principle with speed $n$. 
\end{thm}
This question is still open when $pn\gg 1$. When $p$ is of order one, we should expect to have a large deviation with speed $n^2$  according to the concentraion of measure, but the rate function  should bot be equal to the Gaussian one even when $p=1/2$ because the Dirac at the origin should have rate function bounded above by $\ln p$ (whereas it is infinite in the Gaussian case). 
\subsection{Large deviations for triangle counts}
The traces of Bernoulli matrices have a combinatorial interpretation. For instance, $\tr(\bB_{n}^{3})$ is the number $T_{{n,p}}$ of triangles in the Erd\"os-R\'enyii graph. Observe that its expectation is of order $p^3n^3$. 
In the well known paper \cite[Theorem 4.1]{ChaVa}, it was proved that
\begin{thm} Let 
$$I_{p}(f)=\sup_{\phi}\left\{\int_{0}^{1}\int_{0}^{1} f(x,y)\phi(x,y) dx dy-\frac{1}{2}\int\int \log\left(pe^{2\phi(x,y)}+(1-p) \right) dxdy\right\}$$
and set $\varphi(p,t)=\in\{ I_{p}(f), \int f(x,y)f(y,v)f(v,x)dxdydv\ge 6 t\}$. Then for each $p\in (0,1)$
$$\lim_{n\rightarrow\infty}\frac{1}{n^{2}}\log \Pp\left(T_{n,p}\ge tn^{3}\right) =-\varphi(p,t)\,.$$
\end{thm}
This result extends to any moment $\tr(\bB_n^k)$. However, observe that it does not tell us about deviations of the empirical measure since $x\ra  x^k$ is unbounded so that deviations of the extreme eigenvalues matter. 
It is natural to wonder what happens as well when $p$ goes to zero. This question was attacked in \cite{ChaDe,CoDe,Eldan}, but we state here \cite[Proposition 1.19]{Aug}
\begin{thm} 
 Let $p$ going to zero with $n$ so that $(\log n)^{4}\ll n p^{2}$.  set $v_{n}=n^{2}p^{2}\log(1/p)$. Then for $t\ge 1$
 $$\lim_{n\ra\infty }\frac{1}{v_{n}}\log \Pp\left( \tr( \bB_{n}^{d})\ge t n^{d} p^{d}\right) =-\Phi(t)$$
 where $\Phi(t)=\frac{1}{2}(t-1)^{2/d}$ if $n^{-1}\ll p\ll n^{{-1/2}}$ but $\Phi(t)=
  \min\{ \theta_{t}, \frac{1}{2}(t-1)^{2/d}\}$ if $p\gg n^{{-1/2}}$ and $\theta_{t}$ is the solution of $P_{C_{d}}(\theta_{t})=t$ where $P_{{C_{d}}}$ is the independence polynomial of the $d$-cycle. 
 \end{thm}
\subsection{The Singularity Probability}
A well known problem has been to estimate the probability that a matrix $\tilde \bB_n$ with all independent Bernoulli entries (hence not self-adjoint) is singular. In a breakthrough paper, Tikhomirov \cite{Tikho}, see also \cite{Tikho2}, could exactly estimate it, by showing that the best strategy to achieve singularity is to have a zero column or row vector.
\begin{thm} There exists a finite constant $C$ such that if $C\ln n/n\le p\le \frac{1}{2}$,
$$\mathbb P\left(\tilde\bB_n\mbox{ is singular}\right)= (2+o_n(1))(1-p)^nn$$
\end{thm}
Such an optimal  estimate is not yet known for the symmetric Bernoulli matrix $\bB_n$ (even though it is conjectured) but \cite{CMMM} proves that the probability that it is singular is bounded above by $e^{-O(\sqrt{n})}$. This was improved in an exponential upper bound in \cite{CJMS}.
\section{Open problems}
\begin{enumerate}
\item Local law for Bernoulli matrices when $pn$ is of order one. This could be at best on the scale $\sqrt{n}$ but is tricky even to state  because of the atoms of the limit law.
\item Localization/delocalization of the eigenvectors of Bernoulli matrices for $pn$ of order one (one would conjecture that Dirac masses yield localization but the continuous part yield delocalization, however the right criteria to express this remains to be given).  Find a critical $c^*$ such that for $np>c^*$ there exists delocalized vectors with connected support with high probability.   
\item Large deviations for the empirical measure of the eigenvalues of Bernoulli matrices (all $p$ so that $pn\gg 1$). Even when $p=1/2$ one does not expect to retrieve the Gaussian rate function since the entropy should be finite at $\delta_0$ (as can be seen by requiring all entries to be equal).
\item Precise estimate on the singularity probability in the symmetric case. 
\item In comparison, $d$-regular graphs which are picked uniformly at random are conjectured to be in the universality class of Gaussian random matrices for all $d\ge 3$. This was proven for $d$ going to infinity fast enough \cite{BHKY1,BHKY2}, and recently Huang and Yau could get the local law and the delocalization of the eigenvectors up to $d=3$.
\end{enumerate}
 \bibliographystyle{amsplain}

\bibliography{bibAbel3}

\providecommand{\bysame}{\leavevmode\hbox to3em{\hrulefill}\thinspace}
\providecommand{\MR}{\relax\ifhmode\unskip\space\fi MR }
\providecommand{\MRhref}[2]{%
  \href{http://www.ams.org/mathscinet-getitem?mr=#1}{#2}
}
\providecommand{\href}[2]{#2}
\begin{thebibliography}{10}

\bibitem{ALM}
Amol Aggarwal, Patrick Lopatto, and Jake Marcinek, \emph{{GOE }statistics for
  {L}\'evy matrices}, Arxiv:2002.09355 (2020).

\bibitem{ALY}
Amol {Aggarwal}, Patrick {Lopatto}, and Jake {Marcinek}, \emph{{Eigenvector
  statistics of L\'evy matrices}}, {Ann. Probab.} \textbf{49} (2021), no.~4,
  1778--1846 (English).

\bibitem{ADK1}
Joannes Alt, Raphael Ducatez, and Antti Knowles, \emph{Extremal eigenvalues of
  critical {E}rd{\H o}s-{R}{\'e}nyi graphs}, arXiv 1905.03243 (2019).

\bibitem{ADK2}
\bysame, \emph{Delocalization transition for critical {E}rd{\H o}s-{R}{\'e}nyi
  graphs}, arxiv 2005.14180 (2020).

\bibitem{ADK3}
\bysame, \emph{Poisson statistics and localization at the spectral edge of
  sparse {E}rd{\H o}s-{R}{\'e}nyi graphs}, arxiv 2106.12519 (2021).

\bibitem{AGZ}
Greg~W. Anderson, Alice Guionnet, and Ofer Zeitouni, \emph{An introduction to
  random matrices}, Cambridge Studies in Advanced Mathematics, vol. 118,
  Cambridge University Press, Cambridge, 2010. \MR{2760897}

\bibitem{AZ08}
Greg~W. {Anderson} and Ofer {Zeitouni}, \emph{{A CLT for regularized sample
  covariance matrices}}, {Ann. Stat.} \textbf{36} (2008), no.~6, 2553--2576
  (English).

\bibitem{ArBo}
A.~Arras and C.~Bordenave, \emph{Existence of absolutely continuous spectrum
  for {G}alton-{W}atson random trees}, arXiv:2105.10177 (2021).

\bibitem{ABP}
Antonio Auffinger, G{\'e}rard Ben~Arous, and Sandrine P{\'e}ch{\'e},
  \emph{Poisson convergence for the largest eigenvalues of heavy tailed random
  matrices}, Ann. Inst. Henri Poincar\'e Probab. Stat. \textbf{45} (2009),
  no.~3, 589--610. \MR{2548495}

\bibitem{fanny}
Fanny Augeri, \emph{Large deviations principle for the largest eigenvalue of
  {W}igner matrices without {G}aussian tails}, Electron. J. Probab. \textbf{21}
  (2016), Paper No. 32, 49. \MR{3492936}

\bibitem{Aug}
\bysame, \emph{Nonlinear large deviation bounds with applications to {W}igner
  matrices and sparse {E}rd{\H o}s-{R}{\'e}nyi graphs}, Ann. Probab.
  \textbf{48} (2020), no.~5, 2404--2448.

\bibitem{AGH}
Fanny Augeri, Alice Guionnet, and Jonathan Husson, \emph{Large deviations for
  the largest eigenvalue of sub-{G}aussian matrices}, Comm. Math. Phys.
  \textbf{383} (2021), no.~2, 997--1050.

\bibitem{BHKY1}
Roland {Bauerschmidt}, Jiaoyang {Huang}, Antti {Knowles}, and Horng-Tzer {Yau},
  \emph{{Bulk eigenvalue statistics for random regular graphs}}, {Ann. Probab.}
  \textbf{45} (2017), no.~6A, 3626--3663 (English).

\bibitem{BHKY2}
\bysame, \emph{{Edge rigidity and universality of random regular graphs of
  intermediate degree}}, {Geom. Funct. Anal.} \textbf{30} (2020), no.~3,
  693--769 (English).

\bibitem{BDG}
G.~Ben~Arous, A.~Dembo, and A.~Guionnet, \emph{Aging of spherical spin
  glasses}, Probab. Theory Related Fields \textbf{120} (2001), no.~1, 1--67.
  \MR{1856194}

\bibitem{BAG}
G.~{Ben Arous} and A.~Guionnet, \emph{Large deviations for {W}igner's law and
  {V}oiculescu's non-commutative entropy}, Probab. Theory Rel. \textbf{108}
  (1997), 517--542.

\bibitem{BAG1}
G{\'e}rard Ben~Arous and Alice Guionnet, \emph{The spectrum of heavy tailed
  random matrices}, Comm. Math. Phys. \textbf{278} (2008), no.~3, 715--751.
  \MR{2373441 (2008j:60015)}

\bibitem{BG08}
\bysame, \emph{The spectrum of heavy tailed random matrices}, Comm. Math. Phys.
  \textbf{278} (2008), no.~3, 715--751.

\bibitem{BCK1}
Florent {Benaych-Georges}, Charles {Bordenave}, and Antti {Knowles},
  \emph{{Largest eigenvalues of sparse inhomogeneous {E}rd\H{o}s-{R}{\'e}nyi
  graphs}}, {Ann. Probab.} \textbf{47} (2019), no.~3, 1653--1676 (English).

\bibitem{BG14}
Florent Benaych-Georges and Alice Guionnet, \emph{Central limit theorem for
  eigenvectors of heavy tailed matrices}, Electron. J. Probab. \textbf{19}
  (2014), no. 54, 27. \MR{3227063}

\bibitem{BGM}
Florent Benaych-Georges, Alice Guionnet, and Camille Male, \emph{Central limit
  theorems for linear statistics of heavy tailed random matrices}, Comm. Math.
  Phys. \textbf{329} (2014), no.~2, 641--686. \MR{3210147}

\bibitem{BBG}
Bhaswar~B. Bhattacharya, Sohom Bhattacharya, and Shirshendu Ganguly,
  \emph{Spectral edge in sparse random graphs: {U}pper and lower tail large
  deviations}, Ann. Probab. \textbf{49} (2021), no.~4, --. \MR{4260469}

\bibitem{BCC}
C.~Bordenave, P.~Caputo, and D.~Chaifa, \emph{Spectrum of hermitian heavy
  tailed random matrices}, Comm. Math. Phys. (To appear).

\bibitem{BordCap}
Charles Bordenave and Pietro Caputo, \emph{A large deviation principle for
  {W}igner matrices without {G}aussian tails}, Ann. Probab. \textbf{42} (2014),
  no.~6, 2454--2496. \MR{3265172}

\bibitem{BoCa15}
\bysame, \emph{Large deviations of empirical neighborhood distribution in
  sparse random graphs}, Probab. Theory Related Fields \textbf{163} (2015),
  no.~1-2, 149--222. \MR{3405616}

\bibitem{BG}
Charles Bordenave and Alice Guionnet, \emph{Localization and delocalization of
  eigenvectors for heavy-tailed random matrices}, Probab. Theory Related Fields
  \textbf{157} (2013), no.~3-4, 885--953. \MR{3129806}

\bibitem{BG2}
\bysame, \emph{Delocalization at small energy for heavy-tailed random
  matrices}, Comm. Math. Phys. \textbf{354} (2017), no.~1, 115--159.
  \MR{3656514}

\bibitem{BoLe}
Charles {Bordenave} and Marc {Lelarge}, \emph{{Resolvent of large random
  graphs}}, {Random Struct. Algorithms} \textbf{37} (2010), no.~3, 332--352
  (English).

\bibitem{BoLeSa}
Charles {Bordenave}, Marc {Lelarge}, and Justin {Salez}, \emph{{The rank of
  diluted random graphs}}, {Ann. Probab.} \textbf{39} (2011), no.~3, 1097--1121
  (English).

\bibitem{BSV}
Charles Bordenave, Arnab Sen, and Balint Virag, \emph{Mean quantum
  percolation}, J. Eur. Math. Soc. (JEMS) \textbf{19} (2017), no.~12,
  3679--3707. \MR{3730511}

\bibitem{BouchaudCizeau}
Jean-Philippe Bouchaud and Pierre Cizeau, \emph{Theory of {L}\'evy matrices},
  Phys. Rev. E \textbf{3} (1994), 1810--1822.

\bibitem{CJMS}
M.~Campos, M.~Jenssen, M.~Michelen, and Sahasrabudhe, \emph{The singularity of
  random symmetric matrices is exponentially small}, arxiv 2105.11384 (2021).

\bibitem{CMMM}
M.~Campos, Mattos L., Morris R., and Morrison N., \emph{On the singularity of
  random symmetric matrices}, Duke Math. \textbf{170} (2021), 881--907.

\bibitem{ChaDe}
Sourav Chatterjee and Amir Dembo, \emph{Nonlinear large deviations}, Adv. Math.
  \textbf{299} (2016), 396--450. \MR{3519474}

\bibitem{ChaVa}
Sourav Chatterjee and S.~R.~S. Varadhan, \emph{The large deviation principle
  for the {E}rdos-{R}enyi random graph}, European J. Combin. \textbf{32}
  (2011), no.~7, 1000--1017. \MR{2825532}

\bibitem{CoDe}
Nicholas Cook and Amir Dembo, \emph{Large deviations of subgraph counts for
  sparse {E}rd{\H o}s-{R}{\'e}nyi graphs}, Adv. Math. \textbf{373} (2020),
  107289, 53. \MR{4130460}

\bibitem{Eldan}
Ronen {Eldan}, \emph{{Gaussian-width gradient complexity, reverse log-Sobolev
  inequalities and nonlinear large deviations}}, {Geom. Funct. Anal.}
  \textbf{28} (2018), no.~6, 1548--1596 (English).

\bibitem{EnMe}
Nathana\"{e}l Enriquez and Laurent M\'{e}nard, \emph{Spectra of large diluted
  but bushy random graphs}, Random Structures Algorithms \textbf{49} (2016),
  no.~1, 160--184. \MR{3521277}

\bibitem{ESY}
L\'aszl\'o Erd\H~os, Benjamin Schlein, and Horng-Tzer Yau, \emph{Wegner
  estimate and level repulsion for {W}igner random matrices}, Int. Math. Res.
  Not. IMRN (2010), no.~3, 436--479. \MR{2587574}

\bibitem{EKYY11}
L{\'a}szl{\'o} Erd{\H{o}}s, Antti Knowles, Horng-Tzer Yau, and Jun Yin,
  \emph{Spectral {S}tatistics of {E}rd{\H o}s-{R}\'enyi {G}raphs {II}:
  {E}igenvalue {S}pacing and the {E}xtreme {E}igenvalues}, Comm. Math. Phys.
  \textbf{314} (2012), no.~3, 587--640. \MR{2964770}

\bibitem{EKYY}
\bysame, \emph{The local semicircle law for a general class of random
  matrices}, Electron. J. Probab. \textbf{18} (2013), no. 59, 58. \MR{3068390}

\bibitem{ESY10}
L{\'a}szl{\'o} Erd{\H{o}}s, Benjamin Schlein, and Horng-Tzer Yau, \emph{Wegner
  estimate and level repulsion for {W}igner random matrices}, Int. Math. Res.
  Not. IMRN (2010), no.~3, 436--479.

\bibitem{GZconc}
A.~Guionnet and O.~Zeitouni, \emph{Concentration of the spectral measure for
  large matrices}, Electron. Comm. Probab. \textbf{5} (2000), 119--136.
  \MR{1781846}

\bibitem{GZ3}
A.~Guionnet and O.~Zeitouni, \emph{Large deviations asymptotics for spherical
  integrals}, J. Funct. Anal. \textbf{188} (2002), 461--515.

\bibitem{HuGu}
Alice Guionnet and Jonathan Husson, \emph{Large deviations for the largest
  eigenvalue of {R}ademacher matrices}, Ann. Probab. \textbf{48} (2020), no.~3,
  1436--1465. \MR{4112720}

\bibitem{HeYu}
Yukun He, \emph{Bulk eigenvalue fluctuations of sparse random matrices}, Ann.
  Appl. Probab. \textbf{30} (2020), no.~6, 2846--2879. \MR{4187130}

\bibitem{KY}
Yukun He and Antti Knowles, \emph{Fluctuations of extreme eigenvalues of sparse
  {E}rd{\"o}s-{R}{\'e}nyi graphs}, arXiv 2005.02254 (2020).

\bibitem{HLY}
Jiaoyang {Huang}, Benjamin {Landon}, and Horng-Tzer {Yau}, \emph{{Transition
  from Tracy-Widom to Gaussian fluctuations of extremal eigenvalues of sparse
  {E}rd{\H o}s-{R}{\'e}nyi graphs}}, {Ann. Probab.} \textbf{48} (2020), no.~2,
  916--962 (English).

\bibitem{johansson}
K.~Johansson, \emph{On fluctuations of eigenvalues of random {H}ermitian
  matrices}, Duke Math. J. \textbf{91} (1998), 151--204.

\bibitem{jonsson}
D.~Jonsson, \emph{Some limit theorems for the eigenvalues of a sample
  covariance matrix}, J. Multivariate Anal. \textbf{12} (1982), 1--38.

\bibitem{KSV}
O.~{Khorunzhy}, M.~{Shcherbina}, and V.~{Vengerovsky}, \emph{{Eigenvalue
  distribution of large weighted random graphs}}, {J. Math. Phys.} \textbf{45}
  (2004), no.~4, 1648--1672 (English).

\bibitem{LS}
Ji~Oon {Lee} and Kevin {Schnelli}, \emph{{Local law and Tracy-Widom limit for
  sparse random matrices}}, {Probab. Theory Relat. Fields} \textbf{171} (2018),
  no.~1-2, 543--616 (English).

\bibitem{Tikho2}
Alexander Litvak and Konstantin {Tikhomirov}, \emph{{Singularity of sparse
  random Bernoulli matrices}},  (2020).

\bibitem{ME}
M.~L. Mehta, \emph{Random matrices}, third ed., Pure and Applied Mathematics
  (Amsterdam), vol. 142, Elsevier/Academic Press, Amsterdam, 2004. \MR{2129906
  (2006b:82001)}

\bibitem{PaSh}
Leonid {Pastur} and Mariya {Shcherbina}, \emph{{Eigenvalue distribution of
  large random matrices}}, vol. 171, Providence, RI: American Mathematical
  Society (AMS), 2011 (English).

\bibitem{RuVe}
Mark {Rudelson} and Roman {Vershynin}, \emph{{No-gaps delocalization for
  general random matrices}}, {Geom. Funct. Anal.} \textbf{26} (2016), no.~6,
  1716--1776 (English).

\bibitem{Salez}
Justin {Salez}, \emph{{Every totally real algebraic integer is a tree
  eigenvalue}}, {J. Comb. Theory, Ser. B} \textbf{111} (2015), 249--256
  (English).

\bibitem{ShTi}
Mariya Shcherbina and Brunello Tirozzi, \emph{Central limit theorem for
  fluctuations of linear eigenvalue statistics of large random graphs}, J.
  Math. Phys. \textbf{51} (2010), no.~2, 023523, 20. \MR{2605074 (2011b:05242)}

\bibitem{soshnikov2}
A.~Soshnikov, \emph{Universality at the edge of the spectrum in {W}igner random
  matrices}, Commun. Math. Phys. \textbf{207} (1999), 697--733.

\bibitem{Tala}
M.~Talagrand, \emph{Concentration of measure and isoperimetric inequalities in
  product spaces}, Publ. Math. I.H.E.S. \textbf{81} (1995), 75--203.

\bibitem{TV}
Terence Tao and Van Vu, \emph{The {W}igner-{D}yson-{M}ehta bulk universality
  conjecture for {W}igner matrices}, Electron. J. Probab. \textbf{16} (2011),
  no. 77, 2104--2121. \MR{2851058}

\bibitem{TBT}
E.~{Tarquini}, G.~{Biroli}, and M.~{Tarzia}, \emph{{Level statistics and
  localization transitions of {L}\'evy matrices}}, {Phys. Rev. Lett.}
  \textbf{116} (2016), no.~1, 5 (English), Id/No 010601.

\bibitem{Tikho}
Konstantin {Tikhomirov}, \emph{{Singularity of random Bernoulli matrices}},
  {Ann. Math. (2)} \textbf{191} (2020), no.~2, 593--634 (English).

\bibitem{TiYo}
Konstantin Tikhomirov and Pierre Youssef, \emph{Outliers in spectrum of sparse
  {W}igner matrices}, Random Structures Algorithms \textbf{58} (2021), no.~3,
  517--605. \MR{4234995}

\bibitem{TW5}
C.~A. Tracy and H.~Widom, \emph{Introduction to random matrices}, Lecture Notes
  in Physics, vol. 424, pp.~103--130, Springer, New York, NY, 1993.

\bibitem{TW1}
\bysame, \emph{On orthogonal and symplectic matrix ensembles}, Commun. Math.
  Phys. \textbf{177} (1996), 727--754.

\bibitem{VaVu}
Van Vu, \emph{Recent progress in combinatorial random matrix theory}, arXiv
  2005.02797 (2020).

\bibitem{Wig58}
E.~P. Wigner, \emph{On the distribution of the roots of certain symmetric
  matrices}, Annals Math. \textbf{67} (1958), 325--327.

\bibitem{ZAK}
I.~Zakharevich, \emph{A generalization of {W}igner's law}, Comm. Math. Phys.
  \textbf{268} (2006), 403--414.

\end{thebibliography}

\end{document}